\documentclass[mathpazo]{cicp}

\usepackage{bm}
\usepackage{mathdots}
\usepackage{extarrows}
\allowdisplaybreaks
\usepackage[top=1in,bottom=1in,left=1.25in,right=1.25in]{geometry}
\usepackage{amsfonts}
\usepackage{amsthm}
\usepackage{amssymb,indentfirst}
\usepackage{amsmath}
\usepackage[divps]{xcolor}

\usepackage{mathrsfs}
\usepackage{cases}
\usepackage{titletoc}
\usepackage{multicol}

\usepackage{amssymb}
\def\bc{\begin{center}}       \def\ec{\end{center}}
\def\be{\begin{equation}}     \def\ee{\end{equation}}
\def\ba{\begin{array}}        \def\ea{\end{array}}
\def\bea{\begin{eqnarray}}    \def\eea{\end{eqnarray}}
\def\beaa{\begin{eqnarray*}}  \def\eeaa{\end{eqnarray*}}

\def\ifl{\iffalse}

\begin{document}
%%%%% title : short title may not be used but TITLE is required.
% \title{TITLE}
% \title[short title]{TITLE}
\title[Solving some Navier-Stokes Equations with the initial conditions being some complex-valued periodic functions on $\mathbb{R}^3$]{Solving some Navier-Stokes Equations with the initial conditions being some complex-valued periodic functions on $\mathbb{R}^3$}

%%%%% author(s) :
% single author:
% \author[name in running head]{AUTHOR\corrauth}
% [name in running head] is NOT OPTIONAL, it is a MUST.
% Use \corrauth to indicate the corresponding author.
% Use \email to provide email address of author.
% \footnote and \thanks are not used in the heading section.
% Another acknowlegments/support of grants, state in Acknowledgments section
% \section*{Acknowledgments}
%%%%%

\author[T. Zhang, A. Chen, Fan Bai]{Tao Zhang$^{a}$, \quad  Alatancang Chen$^{a,b,}$\corrauth, \quad Fan Bai$^{a}$}

\address{ $^a$School of Mathematical Sciences of Inner Mongolia University, Hohhot, 010021, China\\
           $^b$Huhhot University for Nationalities, Hohhot, 010051, China}

\email{\tt zhangtaocx@163.com (T. Zhang), alatanca@imu.edu.cn (A. Chen), bf123.student@sina.com (F. Bai).}

% multiple authors:
% Note the use of \affil and \affilnum to link names and addresses.
% The author for correspondence is marked by \corrauth.
% use \emails to provide email addresses of authors
% e.g. below example has 3 authors, first author is also the corresponding
%      author, author 1 and 3 having the same address.
% \author[Zhang Z R et.~al.]{Zhengru Zhang\affil{1}\comma\corrauth,
%       Author Chan\affil{2}, and Author Zhao\affil{1}}
% \address{\affilnum{1}\ School of Mathematical Sciences,
%          Beijing Normal University,
%          Beijing 100875, P.R. China. \\
%           \affilnum{2}\ Department of Mathematics,
%           Hong Kong Baptist University, Hong Kong SAR}
% \emails{{\tt zhang@email} (Z.~Zhang), {\tt chan@email} (A.~Chan),
%          {\tt zhao@email} (A.~Zhao)}
% \footnote and \thanks are not used in the heading section.
% Another acknowlegments/support of grants, state in Acknowledgments section

%%%%% Begin Abstract %%%%%%%%%%%

\begin{abstract}
In this paper, we utilize some series and an iterative method to solve some Navier-Stokes equations with the initial conditions being some complex-valued periodic functions on $\mathbb{R}^3$. Then a new strategy for dealing with the conjecture of the Navier-Stokes equation is given.
\end{abstract}
%%%%% end %%%%%%%%%%%

%%%%% AMS/PACs/Keywords %%%%%%%%%%%
%\pac{...}

%\ams{34L15, 49R05, 34L40, 49J15, 34B05}


%%%% maketitle %%%%%
\maketitle

%%%\tableofcontents

%%%% Start %%%%%%
\section{Introduction}\label{intr}

 Notation

     \begin{equation*}
      \begin{array}{ll}
      \mathbb{R}-\text{the  real  numbers}.\\
       \mathbb{C}-\text{the  complex  numbers}.\\
       \mathbb{R}^{n}=\{(r_1,\cdots,r_n)\mid r_j\in \mathbb{R},\ j=1,2,\cdots,n\}.\\
       \mathbb{N}^{n}=\{(k_{1},\cdots,k_n)\mid k_{j}=0,1,2,\cdots,\ j=1,2,\cdots,n\}.\\
       \mathbb{N}_+^{n}=\{(k_{1},\cdots,k_n)\mid k_{j}=1,2,\cdots,\ j=1,2,\cdots,n\}.\\
       \mathbb{Z}^{n}=\{(k_{1},\cdots,k_{n})\mid \pm k_{j}\in \mathbb{N},\ j=1,2,\cdots,n\}.\\
       %|k|=|k_1+\cdots+k_n|,\ \ \ \ k=(k_1,\cdots,k_n)\in \mathbb{N}^n.\\
       e_1=(2\pi,0,0),\ \ e_2=(0,2\pi,0),\ \ e_3=(0,0,2\pi).\\
       \varphi_{k}=\exp(ik_1 x_1+ik_2 x_2+ik_3x_3),\ \ \ \ k=(k_1,k_2,k_3)\in \mathbb{Z}^3.\\
       \bigwedge_{a,b,c}=\{(k_1,k_2,k_3)\mid(ak_1,bk_2,ck_3)\in \mathbb{N}^3\},\ \ \ \ a,b,c=\pm 1.\\
      % 0^0=1.\\
       %\Omega=\mathbb{R}^3\oplus[0,+\infty).
      \end{array}
      \end{equation*}

    %  \textbf{Definition 1.1.} The space
    %  \begin{equation*}
    %    C^{\infty}(U_n),\ \ \ \ U_n\subseteq \mathbb{R}^n
   %   \end{equation*}
    %   consists of all complex-valued functions $f:\ U_n\rightarrow \mathbb{C}$ such that for each $\beta\in \mathbb{N}^n$,
   %    $D^\beta f$ exists and is continuous on $U_n$.
   %    Then $C^{\alpha}(U_n)$ is a linear space over the field of complex numbers.\vskip8pt

        \textbf{Definition 1.1.}  The set
      \begin{equation*}
    \{f\in C^{\infty}(\mathbb{R}^3)\mid D^{\beta}f=\sum_{k\in \mathbb{Z}^{3}} h_{k}D^{\beta}
    \varphi_{k}, \ \beta\in \mathbb{N}^3,\
   \{h_k\}_{k\in \mathbb{Z}^{3}}\subseteq \mathbb{C}\}
      \end{equation*}
      is a linear space, we use $IE(\mathbb{R}^3)$ to denote this space.\vskip8pt

    The existence and smoothness of the Navier-Stokes equation is an open problem \cite{1}.
Consider the following question with respect to Navier-Stokes equation:\vskip8pt

\textbf{Question 1.2.} Whether there exist $p(x,t), u_j(x,t)\ (u_j(x+e_i,t)=u_j(x,t),$ $i,j=1,2,3)$ such that the following
Navier-Stokes equation hold?
 \begin{numcases}{}
            u_{jt}+\sum\limits_{m=1}^{3}(u_mu_{jx_m}-\nu u_{jx_mx_m})+p_{x_j}=0,\ \ j=1,2,3,\label{Q.1} \\
         u_{1x_1}+u_{2x_2}+u_{3x_3}=0,\ \ x=(x_1,x_2,x_3)\in\mathbb{R}^3,\ t\geq 0,\label{Q.2} \\
         u_j(x,0)=\sum_{k\in \mathbb{Z}^3}A_{jk}\varphi_k\in IE(\mathbb{R}^3),\ \  j=1,2,3,\label{Q.3}
    \end{numcases}
    where $\nu>0$, $u_j(x,0),j=1,2,3$ are real-valued functions.\vskip8pt

\textbf{Theorem 1.3.} Clearly we have
\begin{itemize}
            \item [{\rm (i)}]
           $  A_{jk}=\frac{1}{8\pi^3}\int\limits_{[0,2\pi]^3}u_j(x,0)\varphi_{-k}\text{d}x_1\text{d}x_2\text{d}x_3,\ \ \ \ k\in \mathbb{Z}^3,\ j=1,2,3.$
   \item [{\rm (ii)}]
            $\sum\limits_{k\in \mathbb{Z}^3}|A_{jk}k_1^{m_1}k_2^{m_2}k_3^{m_3}|^2<\infty,\ \ \ \ j=1,2,3,\ (m_1,m_2,m_3)\in \mathbb{N}^3$.
  \item [{\rm (iii)}]  By the conditions \eqref{Q.2} and \eqref{Q.3} we get
  \begin{equation*}
     \sum\limits_{k\in \mathbb{Z}^3}\sum\limits_{j=1}^3 A_{jk}k_j\varphi_k=0,\ \ \ \ \ \ \ k\in \mathbb{Z}^3.
  \end{equation*}
  Note that the sequence $\{\varphi_k\}_{k\in \mathbb{Z}^3}$ is linearly independent, so we have
  \begin{equation*}
     \sum\limits_{j=1}^3A_{jk}k_j=0,\ \ \ \ \ \ \ k\in \mathbb{Z}^3.
  \end{equation*}
   \item [{\rm (iv)}]
 $\sum\limits_{k\in \bigwedge_{a,b,c}}A_{jk}\varphi_k\in IE(\mathbb{R}^3),\ \ \ \ j=1,2,3,\ a,b,c=\pm 1.$
\end{itemize}\vskip8pt

    In this paper, base on the idea of paper \cite{4}, we can solve the following PDEs in some cases:
    \begin{equation}\label{QQ}
      \left\{
      \begin{array}{l}
          \text{the PDEs \eqref{Q.1} and \eqref{Q.2}},\\
          u_j(x,0)=\sum\limits_{k\in \bigwedge_{a,b,c}}B_{jk}\varphi_{k}\in IE(\mathbb{R}^3),\ \  j=1,2,3,
      \end{array}
      \right.
    \end{equation}
    where $a,b,c=\pm 1$.\vskip8pt

 If we let
 \begin{equation}\label{QE}
  \sum\limits_{k\in \mathbb{Z}^3}A_{jk}\varphi_{k}=\sum\limits_{a,b,c=\pm 1}\ \sum\limits_{k\in \bigwedge_{a,b,c}}B_{jk}\varphi_{k}.
 \end{equation}
 Then there should be some relation between the solution of the PDEs \eqref{Q.1}-\eqref{Q.3} and the solution of the PDEs \eqref{QQ}.
So a new strategy for dealing with the conjecture of the Navier-Stokes equation is given.

%\section{Preliminaries}

 %     \textbf{Definition 2.1.} The space
%%      \begin{equation*}
%      C(U)
%      \end{equation*}
%       consists of all complex-valued continuous functions on $U\subseteq \mathbb{R}^n$.

%%--------------------------------------------------------
\section{Main results}
%%--------------------------------------------------------

   First we solve the following PDEs:
   \begin{equation}\label{sss}
     \left\{
     \begin{array}{l}
       \text{the PDEs \eqref{Q.1} and \eqref{Q.2}},\\
         u_j(x,0)=\sum\limits_{k\in \mathbb{N}^3}B_{jk}\varphi_{k}\in IE(\mathbb{R}^3),\ \  j=1,2,3.
     \end{array}
     \right.
   \end{equation}\vskip8pt

   Suppose that the PDEs \eqref{sss} has a solution satisfying:
    \begin{equation}\label{s.4}
    \left\{
    \begin{array}{l}
      u_j(x,t)=\sum\limits_{k\in \mathbb{N}^3}T_{jk}(t)\varphi_k,\ \ j=1,2,3;\\
      p(x,t)=\sum\limits_{k\in \mathbb{N}^3}T_{4k}(t)\varphi_k.
    \end{array}
    \right.
    \end{equation}
    Assume that the series \eqref{s.4} satisfies the following conditions:
     \begin{numcases}{}
         u_j=\sum\limits_{k\in \mathbb{N}^3}T_{jk}(t)\varphi_k\in C(\mathbb{R}^3\oplus[0,+\infty)),\ \ \ \ j=1,2,3, \label{11nb}\\
        p=\sum\limits_{k\in \mathbb{N}^3}T_{4k}(t)\varphi_{k}\in C(\mathbb{R}^3\oplus[0,+\infty)), \label{22nb}\\
        u_{jt}=\sum\limits_{k\in \mathbb{N}^3}T'_{jk}(t)\varphi_{k}\in C(\mathbb{R}^3\oplus[0,+\infty)),\  \ \ \ j=1,2,3, \label{33nb}\\
        u_{jx_m}=\sum\limits_{k\in \mathbb{N}^3} ik_mT_{jk}(t)\varphi_{k}\in C(\mathbb{R}^3\oplus[0,+\infty)),\  \ \ \ m,j=1,2,3, \label{44nb}\\
         u_{jx_mx_m}=\sum\limits_{k\in \mathbb{N}^3} -k_m^2T_{jk}(t)\varphi_{k}\in C(\mathbb{R}^3\oplus[0,+\infty)),\  \ \ \ m,j=1,2,3, \label{55nb}\\
         p_{x_j}=\sum\limits_{k\in \mathbb{N}^3}ik_jT_{4k}(t)\varphi_{k}\in C(\mathbb{R}^3\oplus[0,+\infty)),\  \ \ \ j=1,2,3, \label{66nb}\\
         u_mu_{jx_m}=\sum\limits_{k\in \mathbb{N}^3}\eta_{mjk}\varphi_{k} \in C(\mathbb{R}^3\oplus[0,+\infty)),\ \ \ \  m,j=1,2,3, \label{77nb}
    \end{numcases}
    where
    \begin{equation*}
    %\begin{array}{r}
     \eta_{mjk}=\sum\limits_{k^{[1]}+k^{[2]}=k}k^{[2]}_{m}
     T_{mk^{[1]}}T_{jk^{[2]}},\ \
      k^{[l]}=(k_{1}^{[l]},k_{2}^{[l]},k_{3}^{[l]})\in \mathbb{N}^3,\ \
      l=1,2,\ m,j=1,2,3.
    %\end{array}
    \end{equation*}
    Then substituting the series \eqref{s.4} into the equations \eqref{sss} we get
    \begin{equation*}
    \left\{
      \begin{array}{l}
      %\begin{array}{r}
         T'_{m,(0,0,0)}+ \sum\limits_{k> (0,0,0)}[T'_{mk}+
         \sum\limits_{j=1}^{3}(\sum\limits_{k^{[1]}+k^{[2]}=k}ik^{[2]}_{j}T_{jk^{[1]}}T_{mk^{[2]}}
    +\nu k_j^2T_{mk})
    +ik_m T_{4k}]\varphi_{k}=0,\ \
     m=1,2,3,\\
     % \end{array}\\
     \sum\limits_{k\in \mathbb{N}^3}(ik_1 T_{1k}+ik_2 T_{2k}+ik_3 T_{3k})\varphi_{k}=0,\\
     u_j(x,0)=\sum\limits_{k\in \mathbb{N}^3}B_{jk}\varphi_{k}=\sum\limits_{k\in \mathbb{N}^3}T_{jk}(0)\varphi_k,\ \  j=1,2,3.
      \end{array}
      \right.
    \end{equation*}
    Note that the sequence $\{\varphi_{k}\}_{k\in \mathbb{N}^3}$ is linearly independent, so the above equations are equivalent to
    the following ODEs:
    \begin{equation*}
      \left\{
      \begin{array}{c}
          T'_{j,(0,0,0)}=0,\\
          T_{j,(0,0,0)}(0)=B_{j,(0,0,0)},
      \end{array}
      \right.\ \   j=1,2,3,
    \end{equation*}
    and
    \begin{equation*}
      \left\{
      \begin{array}{c}
         T'_{mk}+ \sum\limits_{j=1}^{3}(\sum\limits_{k^{[1]}+k^{[2]}=k}ik^{[2]}_{j}T_{jk^{[1]}}T_{mk^{[2]}}
         +\nu k_j^2T_{mk})+ik_m T_{4k}=0,\ \ m=1,2,3,\\
         k_1 T_{1k}+k_2 T_{2k}+k_3 T_{3k}=0,\\
          T_{jk}(0)=B_{jk},\ \ \ \ \ \ \ \ j=1,2,3,
      \end{array}
      \right.
    \end{equation*}
    where $k>(0,0,0)$.
    By the equations $k_1 T_{1k}+k_2 T_{2k}+k_3 T_{3k}=0$, $k>(0,0,0)$, we have
    \begin{equation*}
         k_1T'_{1k}+k_2 T'_{2k}+k_3 T'_{3k}=0,\ \ \ \ \ \  k>(0,0,0).
    \end{equation*}
    Then we obtain
    \begin{equation*}
    %\begin{array}{r}
      \sum\limits_{m=1}^{3}k_{m}\sum\limits_{j=1}^{3}\sum\limits_{k^{[1]}+k^{[2]}=k,\atop k^{[1]},k^{[2]}>(0,0,0)}
     ik^{[2]}_{j}T_{jk^{[1]}}T_{mk^{[2]}}+T_{4k}\sum\limits_{m=1}^{3}ik_m^2
    =0,\ \  k>(0,0,0).
    %\end{array}
    \end{equation*}
   So we get
   \begin{equation*}
    \left\{
    \begin{array}{l}
      T_{j,(0,0,0)}(t)=B_{j,(0,0,0)},\ \ \ \ \ \ \ \ j=1,2,3,\\
     T_{4,(0,0,0)}=a,\ \ \ \ \ \ \ \  a\ \text{is an arbitrary constant,}\\
    T_{4k}(t)=\frac{-\sum\limits_{m=1}^{3}k_m\sum\limits_{j=1}^{3}\ \sum\limits_{k^{[1]}+k^{[2]}=k,\ k^{[1]},k^{[2]}>(0,0,0)}
     k^{[2]}_{j}T_{jk^{[1]}}T_{mk^{[2]}}}{\sum\limits_{m=1}^{3}k_m^2}, \ \ \ \ k>(0,0,0),\\
    % \begin{array}{l}
    T_{jk}(t)=\exp(-P_{k}(t))(\int_{0}^{t}Q_{jk}(s)\exp(P_{k}(s))\text{d}s+B_{jk}),\ \
   j=1,2,3, \ \ k>(0,0,0),
  %   \end{array}
    \end{array}
    \right.
    \end{equation*}
    where
    \begin{equation*}
      \left\{
      \begin{array}{l}
      %\begin{array}{r}
      Q_{mk}=- \sum\limits_{j=1}^{3}\ \sum\limits_{k^{[1]}+k^{[2]}=k,\ k^{[1]},k^{[2]}>(0,0,0)}
     ik^{[2]}_{j}T_{jk^{[1]}}T_{mk^{[2]}}-ik_m T_{4k},\ \ \ \ m=1,2,3,\\
     % \end{array}
     P_{k}(t)=\sum\limits_{j=1}^{3}ik_jB_{j,(0,0,0)}t+\nu \sum\limits_{j=1}^{3}k_j^2t.
      \end{array}
      \right.
    \end{equation*}\vskip8pt

     Clearly we have:\vskip8pt

     \textbf{Theorem 2.1.}   If the series \eqref{s.4} we obtain satisfies the conditions
     \eqref{11nb}-\eqref{77nb}, then it is a solution of the PDEs  \eqref{sss}.\vskip8pt

     By the Abel identities \cite{ci} we have:\vskip8pt

      \textbf{Lemma 2.2.} For any $k=1,2,\cdots$, we have
      \begin{itemize}
            \item [{\rm (i)}]
       $\sum\limits_{m=1}^{k}\frac{(k+1)!}{m!(k+1-m)!}
                            m^m(k+1-m)^{k-m}=k(k+1)^k,$

            \item [{\rm (ii)}] $\sum\limits_{m=1}^{k}\frac{(k+1)!}{m!(k+1-m)!}
                            m^{m-1}(k+1-m)^{k-m}=2k(k+1)^{k-1}\leq 2(k+1)^{k}.$
             \end{itemize}\vskip8pt

       \textbf{Corollary 2.3.} For any $ k=(k_1,\cdots,k_n)\in \mathbb{N}_+^n$, we have
      \begin{equation*}
\sum\limits_{(1,\cdots,1)\leq (m_1,\cdots,m_n)\leq k}\ m_1\prod\limits_{j=1}^n \frac{m_j^{m_j-1}(k_j+1-m_j)^{k_j-m_j}}{m_j!(k_j+1-m_j)!}\leq 2^{n-1}k_1\prod\limits_{j=1}^n\frac{(k_j+1)^{k_j}}{(k_j+1)!}.
      \end{equation*}\vskip8pt

     %  \textbf{Theorem 3.4.} Let $|k|=|k_1+k_2+k_3|$. Then there exists $n=(n_1,n_2,n_3)\in \mathbb{N}^3$ such that
    % \begin{equation}\label{guina0}
   %  |T_{ik}(t)|\leq \prod\limits_{j=1}^3\frac{k_j^{k_j-1}}{k_j!}e^{-\nu|k|t},\ \ \ \ k\geq n,\  i=1,2,3.
   %  \end{equation}\vskip8pt

    %  \textbf{Proof.} By Theorem 1.3. (ii), there exists some $k_0\in \mathbb{N}^3,\ k_0>(3,3,3)$ such that
    %  \begin{equation*}
    %    |A_{jk}|\leq \frac{1}{|k|^2},\ \ \ \ \ \ \ k>k_0,\ j=1,2,3.
    %  \end{equation*}

     \textbf{Lemma 2.4.} Let $k=(k_1,k_2,k_3)\in \mathbb{N}^3,$ $|k|=|k_1+k_2+k_3|,\ \nu\geq 1$. If
       \begin{equation*}
         |B_{jk}|\leq  \frac{e^{-|k|}}{10^3}\prod\limits_{j=1,2,3,\ k_j>0}\frac{k_j^{k_j-1}}{k_j!},\ \ \ \ \ \ \ k>(0,0,0),\ j=1,2,3.
       \end{equation*}
      Then we have
     \begin{equation}\label{guina1}
     |T_{ik}(t)|\leq \frac{1}{100}\prod\limits_{j=1,2,3,\ k_j>0}\frac{k_j^{k_j-1}}{k_j!}\exp(-\nu|k|t-|k|),\ \ \ \ k>(0,0,0),\  i=1,2,3.
     \end{equation}\vskip8pt

    \textbf{Proof.} We prove the inequalities \eqref{guina1} by the induction method.
    By a simple calculate we can induce that the inequalities \eqref{guina1} hold when $|k|=1,2$.
    Suppose that it hold for any $|k|<n$ $(n>2)$, then by Corollary 2.3, for any $k=(k_1,k_2,k_3)\geq(1,1,1)$, $|k|=n$ (without loss of generality we suppose that $k_1\geq k_2,k_3$), we have
    \begin{equation*}
     % \left\{
      \begin{array}{l}
        |T_{4k}(t)|\leq\frac{\sum\limits_{m=1}^{3}k_m\sum\limits_{j=1}^{3}\ \sum\limits_{k^{[1]}+k^{[2]}=k,\ k^{[1]},k^{[2]}>(0,0,0)}
     k^{[2]}_{j}|T_{jk^{[1]}}||T_{mk^{[2]}}|}{\sum\limits_{m=1}^{3}k_m^2}\\
     \ \ \ \  \ \ \ \ \ \ \ \ \ \leq\frac{3\sum\limits_{m=1}^{3}k_m}{10^4\sum\limits_{m=1}^{3}k_m^2}(
     \sum\limits_{1\leq m_i\leq k_i,i=1,2,3}\ m_1\prod\limits_{j=1}^3 \frac{m_j^{m_j-1}(k_j-m_j)^{k_j-m_j-1}}{m_j!(k_j-m_j)!}+ 6 k_1 \prod\limits_{j=1}^3\frac{k_j^{k_j-1}}{k_j!})\exp(-\nu|k|t-|k|)\\
     \ \ \ \  \ \ \ \ \ \ \ \ \ \leq\frac{60}{10^4}\ \prod\limits_{j=1}^3\frac{k_j^{k_j-1}}{k_j!}\exp(-\nu|k|t-|k|),\\
     % \ \ \ \  \ \ \ \ \ \ \ \ \ \leq\frac{21}{10^4} \prod\limits_{j=1}^3\frac{k_j^{k_j-1}}{k_j!}\exp(-\nu|k|t-|k|),\\
         |Q_{jk}(s)|\leq\sum\limits_{j=1}^{3}\ \sum\limits_{k^{[1]}+k^{[2]}=k,\ k^{[1]},k^{[2]}>(0,0,0)}
     k^{[2]}_{j}|T_{jk^{[1]}}(s)||T_{mk^{[2]}}(s)|+k_j |T_{4k}(s)|\\
      \ \ \ \  \ \ \ \ \ \ \ \ \ \leq\frac{90k_1}{10^4}\prod\limits_{j=1}^3\frac{k_j^{k_j-1}}{k_j!}\exp(-\nu|k|t-|k|),\ \ \ \ j=1,2,3,\\
    |T_{jk}(t)|\leq\exp(-\nu \sum\limits_{j=1}^{3}k_j^2t)(\int_{0}^{t}|Q_{jk}(s)|\exp(\nu \sum\limits_{j=1}^{3}k_j^2s)\text{d}s+|B_{jk}|)\\
     \ \ \ \  \ \ \ \ \ \ \ \ \ \leq\frac{1}{100}\prod\limits_{j=1}^3\frac{k_j^{k_j-1}}{k_j!}\exp(-\nu|k|t-|k|),\ \ \ \ j=1,2,3.
     \end{array}
     \end{equation*}
    In a similar way, we can prove that the inequalities \eqref{guina1} hold for any $k=(k_1,k_2,k_3)\in \mathbb{N}^3$ with $k_1=0$ or $k_2=0$ or $k_3=0$.
\vskip8pt

  \textbf{Theorem 2.5.} Let $|k|=|k_1+k_2+k_3|,\ \nu\geq 1$. If
       \begin{equation*}
         |B_{jk}|\leq  \frac{e^{-|k|}}{10^3}\prod\limits_{j=1,2,3,\ k_j>0}\frac{k_j^{k_j-1}}{k_j!},\ \ \ \ \ \ \ k>(0,0,0),\ j=1,2,3.
       \end{equation*}
     Then the series \eqref{s.4} we obtain is a solution of the PDEs \eqref{sss}.\vskip8pt

\textbf{Proof.} We only need to prove that the series \eqref{s.4} we obtain satisfies the conditions \eqref{11nb}- \eqref{77nb}.
 Note that $\frac{k^m}{m!}\leq e^{k},\ m\in \mathbb{N}$, so we have
\begin{equation*}
\left\{
\begin{array}{l}
  |T_{ik}(t)\varphi_k|\leq \frac{1}{100}e^{-\nu|k|t},\ \ \ \ k>(0,0,0),\  i=1,2,3,\\
  |T_{4k}(t)\varphi_k|\leq \frac{60}{10^4}e^{-\nu|k|t},\ \ \ \ k>(0,0,0).
\end{array}
\right.
\end{equation*}
Hence the series \eqref{s.4} converges absolutely on $\Omega\oplus[0,+\infty)$.
    It means that the series \eqref{s.4} we obtain satisfies the conditions \eqref{11nb}-\eqref{22nb}. Moreover, we can prove that
    \begin{equation*}
    \left\{
    \begin{array}{l}
       |T'_{mk}\varphi_k|=|\sum\limits_{j=1}^{3}(\sum\limits_{k^{[1]}+k^{[2]}=k}ik^{[2]}_{j}T_{jk^{[1]}}T_{mk^{[2]}}
         +\nu k_j^2T_{mk})+ik_m T_{4k}||\varphi_k|\\
     \ \ \ \ \ \ \ \ \ \ \ \ \  \leq (\nu |k|^2|T_{mk}|+|Q_{mk}|)<(\frac{\nu |k|^2}{100}+\frac{90|k|}{10^4})e^{-\nu|k|t},\ \ \ \ m=1,2,3,\ k>(0,0,0),\\
     |k_j^2 T_{mk}\varphi_k|\leq \frac{|k|^2}{100}e^{-\nu|k|t},\ \ \ \ \ \ \ \ \ \ \ \ m,j=1,2,3,\ k>(0,0,0),\\
     |k_j T_{4k}\varphi_k|\leq \frac{60|k|}{10^4}e^{-\nu|k|t},\ \ \ \ \ \ \ \ \ \ \ \ \ j=1,2,3,\ k>(0,0,0),\\
     |\eta_{jmk}\varphi_k|\leq \frac{10}{10^4}e^{-\nu|k|t},\ \ \ \ \ \ \ \ \ \ \ \ \ \ \  m,j=1,2,3,\ k>(0,0,0).
    \end{array}
    \right.
    \end{equation*}
    So the series \eqref{s.4} we obtain satisfies the conditions \eqref{33nb}-\eqref{77nb}.
    Therefore it is a solution of the equations \eqref{sss} by Theorem 2.1.\vskip8pt

    Similarly, we can get:\vskip8pt

   \textbf{Theorem 2.6.}  Let $\nu\geq 1$. For any $a,b,c=\pm 1$, if
       \begin{equation*}
         |B_{jk}|\leq  \frac{\exp(-\sum_{j=1}^3|k_j|)}{10^3}\prod\limits_{j=1,2,3,\ |k_j|>0}\frac{|k_j|^{|k_j|-1}}{|k_j|!},\ \ \ \ \  (ak_1,bk_2,ck_3)>0,\ j=1,2,3,
       \end{equation*}
    then the solution of the PDEs \eqref{QQ} exists.\vskip8pt

\textbf{Theorem 2.7.} For any $a,b,c=\pm 1$, if the functions $p$, $u_j,j=1,2,3$ satisfy the PDEs  \eqref{QQ},
then $\overline{p}$, $\overline{u_j},j=1,2,3$ satisfy the following PDEs
\begin{equation*}
  \left\{
  \begin{array}{l}
     \text{the PDEs \eqref{Q.1} and \eqref{Q.2}},\\
      u_j(x,0)=\sum\limits_{k\in \bigwedge_{a,b,c}}\overline{B_{jk}}\overline{\varphi_{k}},\ \  j=1,2,3.
  \end{array}
  \right.
\end{equation*} \vskip8pt

Next suppose that the equality \eqref{QE} holds, then we give the following conjectures:\vskip8pt

\textbf{Conjecture 2.8.}  If for any $a,b,c=\pm 1$, the solution of the PDEs \eqref{QQ}
exists (and unique), then the solution of the PDEs \eqref{Q.1}-\eqref{Q.3} exists (and unique).\vskip8pt

\textbf{Conjecture 2.9.} Suppose that for any $a,b,c=\pm 1$, the functions $p_{abc}$, $u_{jabc},j=1,2,3$ satisfy the PDEs \eqref{QQ},
  and that the functions $p,u_j,j=1,2,3$ satisfy the PDEs \eqref{Q.1}-\eqref{Q.3}, then there exist some nonlinear functions $T_j,\ j=1,2,3,4,$ such that
   \begin{equation*}
   \left\{
   \begin{array}{l}
      T_j(\Delta)=u_j(x,t),\ \ \ \ j=1,2,3,\\
     T_4(\Delta)=p(x,t),
   \end{array}
   \right.
   \end{equation*}
where $\Delta=\{p_{abc}(x,t), u_{jabc}(x,t)\mid j=1,2,3,\ a,b,c=\pm 1\}$.

\section*{Acknowledgments}

The paper is supported by the Natural Science Foundation of China (no. 11371185) and the Natural Science Foundation of Inner Mongolia,
China (no. 2013ZD01).

%%%% Bibliography  %%%%%%%%%%

%%\vfill \hfill \fbox{\today}
\end{document}